\documentclass{amsart}
\usepackage{amsmath, amsthm, amsfonts, amssymb, txfonts, graphicx, hyperref, verbatim}
\theoremstyle{plain}
\newtheorem{thrm}{Theorem}[section]

\newtheorem{prpstn}[thrm]{Proposition}

\newtheorem*{cnjctr}{Conjecture}
\newtheorem{hypthss}{Hypothesis}

\numberwithin{sblmm}{thrm} 
\numberwithin{equation}{section}

\begin{document}
\title{3-tuples have at most 7 prime factors infinitely often}
\author{James Maynard}
\address{Mathematical Institute, 24–-29 St Giles', Oxford, OX1 kLB}
\email{maynard@math.ox.ac.uk}
\thanks{Supported by EPSRC Doctoral Training Grant EP/P505216/1 }
\date{}
\begin{abstract}
Let $L_1$, $L_2$ $L_3$ be integer linear functions with no fixed prime divisor. We show there are infinitely many $n$ for which the product $L_1(n)L_2(n)L_3(n)$ has at most 7 prime factors, improving a result of Porter. We do this by means of a weighted sieve based upon the Diamond-Halberstam-Richert multidimensional sieve.
\end{abstract}
\maketitle
\section{Introduction}
We consider a set of integer linear functions
\begin{equation}
L_i(x)=a_ix+b_i,\qquad i\in\{1,\dots,k\}.
\end{equation}
We say such a set of functions is \textit{admissible} if their product has no fixed prime divisor. That is, for every prime $p$ there is an integer $n_p$ such that none of $L_i(n_p)$ are a multiple of $p$. We are interested in the following conjecture.
\begin{cnjctr}[Prime $k$-tuples Conjecture]
Given an admissible set of integer linear functions $L_i(x)$ ($i\in\{1,\dots,k\}$), there are infinitely many integers $n$ for which all the $L_i(n)$ are prime.
\end{cnjctr}
With the current technology it appears impossible to prove any case of the prime $k$-tuples conjecture for $k\ge 2$.

Although we cannot prove that the functions are simultaneously prime infinitely often, we are able to show that they are \textit{almost prime} infinitely often, in the sense that their product has only a few prime factors. This was most notably achieved by Chen \cite{Chen} who showed that there are infinitely many primes $p$ for which $p+2$ has at most 2 prime factors. His method naturally generalises to show that for a pair of admissible functions the product $L_1(n)L_2(n)$ has at most 3 prime factors infinitely often.

Similarly sieve methods can prove analogous results for any $k$. We can show that the product of $k$ admissible functions $L_1(n)\dots L_k(n)$ has at most $r_k$ prime factors infinitely often, for some explicitly given value of $r_k$. We see that the prime $k$-tuples conjecture is equivalent to showing we can have $r_k=k$ for all $k$. The current best values of $r_k$ grow asymptotically like $k\log{k}$ and explicitly for small $k$ we can take $r_2=3$ (Chen, \cite{Chen}), $r_3=8$ (Porter, \cite{Porter}), $r_4=11$, $r_5=15$, $r_6=18$, $r_7=22$, $r_8=26$, $r_9=30$, $r_{10}=34$ (Maynard, \cite{Maynard}). %$r_4=12$, $r_5=16$, $r_6=20$ (Diamond and Halberstam \cite{DiamondHalberstam}), $r_7=24$, $r_8=28$, $r_9=33$, $r_{10}=38$ (Ho and Tsang, \cite{HoTsang})} 
Heath-Brown \cite{HeathBrown} showed that infinitely often there are $k$-tuples where all the functions $L_i$ have individually at most $C\log{k}$ prime factors, for an explicit constant $C$.
\section{Main Results}
\begin{thrm}\label{thrm:MainTheorem}
Let $\mathcal{L}=\{L_1,L_2,L_3\}$ be an admissible 3-tuple of integer linear functions. Then there are infinitely many $n$ for which the product $L_1(n)L_2(n)L_3(n)$ has at most 7 prime factors.
\end{thrm}
Results such as Theorem \ref{thrm:MainTheorem} which show that $k$-tuples take few prime factors infinitely often tend to use weighted sieves. When $k=2$ the best results fix one of the functions to attain a prime value, and use a $k-1$ dimensional sieve to show that the remaining functions have few prime factors infinitely often. When $k\ge 4$ better bounds are obtained by instead using a $k$-dimensional sieve on all the factors. When $k=3$ both approaches can show that a 3-tuple has at most $8$ prime factors infinitely often. In this paper we mix these two approaches, using $k$-dimesional sieves to estimate some terms and $k-1$ dimensional sieves to estimate other terms. This is the key innovation which allows us to reduce the bound $r_3$ from $8$ to $7$.

A similar mixed approach can be used to gain an improvement over the traditional approach when $k\ge 4$, but in this case we find it is superior to use a slightly different argument based on the Selberg sieve, and in the author's work \cite{Maynard}.

Theorem \ref{thrm:MainTheorem} follows almost immediately from the following proposition, which may be viewed as a sharpening of Halberstam and Diamond's Theorem 11.1 in \cite{DiamondHalberstam}, tailored to our specific application.
\begin{prpstn}\label{prpstn:MainProp}
Let $\mathcal{L}$ be and admissible $k$-tuple of integer linear functions. Let $r_k$ be a natural number satisfying
\begin{align*}
r_k>\max(N(u,v;k),uk-1)
\end{align*}
where
\begin{align*}
N(u,v;k)&=u k-1+\frac{k}{f_k(v)}\left(I_1-I_2-\frac{e^\gamma(u-1)}{v}f_{k-1}\left(\frac{v}{2}\right)\right),\\
I_1&=\int_{1/v}^{1/u}\min\left(F_k(v-vs),e^\gamma F_{k-1}\left(\frac{v}{2}\right)w(v-vs)\right)\frac{1-us}{s}ds,\\
I_2&=\int_{1/u}^{1-1/v}\max\left(f_k(v-vs),e^\gamma f_{k-1}\left(\frac{v}{2}\right)w(v-vs)\right)\frac{us-1}{s}ds,
\end{align*}
and $u,v$ are real numbers satisfying $v/(v-1)<u<v$ and $v>\beta_k$. 

Then there are infinitely many $n$ for which the product $L_1(n)\dots L_k(n)$ has at most $r_k$ prime factors.

Here $F_k$ and $f_k$ are the upper and lower sieve functions from the $k$-dimensional Diamond-Halberstam-Richert sieve (described in detail in \cite{DiamondHalberstam}[Chapter 6]), and $\beta_k$ is the sifting limit of $f_k$. $w(u)$ is the Buchstab function defined by the delay differential equation given by \eqref{eq:BuchstabFn}.
\end{prpstn}
In \cite{DiamondHalberstam}[Theorem 11.1] one has a similar result but instead of $N(u,v;k)$ one has the expression
\begin{equation}
uk-1+\frac{k}{f_k(v)}\int_{1/v}^{1/u}F_k(v-vs)\frac{1-us}{s}ds.
\end{equation}
We can see that certainly whenever $f_{k-1}(v/2)>0$ Proposition \ref{prpstn:MainProp} gives a superior bound. The optimal choice of $u$ and $v$ when using \cite{DiamondHalberstam}[Theorem 11.1] with $k=3$ is approximately $u=1.5$, $v=12$, and since with these values $f_2(v/2)>0$ we expect a small improvement.

Using numerical integration we can establish Theorem 2.1 from Proposition \ref{prpstn:MainProp}. We consider a 3-tuple of integer linear functions, and take $u=2$, $v=12$. We find that
\begin{equation}
N(2,12;3)=6.943\dots
\end{equation}
and so we may take $r_k=7$ in Proposition \ref{prpstn:MainProp}.

We note that the proof of Proposition \ref{prpstn:MainProp} makes much use of the fact we are dealing with integer linear functions for which much more is known about the distribution of the values they take - Diamond and Halberstam's Theorem 11.1 holds in much more general circumstances. This also gives no change to the asymptotic bound of $k\log{k}+O(k)$ for the number of prime factors of a $k$-tuple when $k$ is large.
\section{Proof of Proposition \ref{prpstn:MainProp}}
To simplify the argument we adopt a normalisation of our functions, as done originally by Heath-Brown \cite{HeathBrown}. Since we are only interested in the showing any admissible $k$-tuple has at most $r_k$ prime factors infinitely often (for some explicit $r_k$), by considering the functions $L_i(An+B)$ instead of $L_i(n)$ for suitably chosen constants $A$ and $B$, we may assume without loss of generality that our functions satisfy the following hypothesis.
\begin{hypthss}\label{hypthss:Normalised}
$\mathcal{L}=\{L_1,\dots,L_k\}$ is an admissible $k$-tuple of linear functions. The functions $L_i(n)=a_i n+b_i$ ($1\le i \le k$) are distinct with $a_i>0$. Each of the coefficients $a_i$ is composed of the same primes none of which divides the $b_j$. If $i\ne j$, then any prime factor of $a_i b_j-a_j b_i$ divides each of the $a_l$.
\end{hypthss}
Let $\mathcal{L}=\{L_1,\dots,L_k\}$ be an admissible $k$-tuple satisfying Hypothesis \ref{hypthss:Normalised}. We view this $k$-tuple as fixed, and so any implied constants from $\ll$ or $O$-notation may depend on the $k$-tuple without explicit reference. We define
\begin{align}
\Pi(n)&=\prod_{i=1}^k L_i(n),\\
A&=\prod_{i=1}^{k}a_i,\\
v_{p}(\mathcal{L})&=\#\{1\le n< p:\Pi(n)\equiv 0\pmod{p}\}.
\end{align}
By Hypothesis \ref{hypthss:Normalised} we have
\begin{equation}
v_p(\mathcal{L})=\begin{cases}
k,\qquad &p\nmid A,\\
0,&p|A.
\end{cases}
\end{equation}
Finally we define the quantity
\begin{equation}
P(z)=\prod_{p<z}p.
\end{equation}
We consider the sum
\begin{align}
S=S(\tau;N,z)=\sum_{\substack{N\le n< 2N\\ (\Pi(n),P(z))=1}}(\tau-\Omega(\Pi(n))).
\end{align}
If for some fixed constant $\tau$ we can show $S(\tau;N,z)>0$ (for suitable $z$, $N$), then there must be an $n\in[N,2N)$ such that $\Omega(\Pi(n))<\tau$. Thus if we can show $S(\tau;z,N)>0$ for all sufficiently large $N$ (with suitable $z$ depending on $N$) then there are infinitely many $n$ such that $\Omega(\Pi(n))\le \lfloor\tau\rfloor$.

We first split the sum $S$ up as a weighted sum over the prime factors of each of the functions $L_j(n)$, based on a new parameter $y$.
\begin{align}
S&=\sum_{\substack{N\le n< 2N\\ (\Pi(n),P(z))=1}}(\tau-\Omega(\Pi(n)))\nonumber\\
&=\sum_{\substack{N\le n< 2N\\ (\Pi(n),P(z))=1}}\left(\tau-\frac{\log{\Pi(n)}}{\log{y}}-\sum_{p|\Pi(n)}\left(1-\frac{\log{p}}{\log{y}}\right)\right)+O\Biggl(\sum_{\substack{N\le n< 2N\\(\Pi(n),P(z))=1\\\Pi(n)\text{ not square-free}}}\log{N}\Biggr)\nonumber\\
&=\sum_{\substack{N\le n< 2N\\ (\Pi(n),P(z))=1}}\left(\tau-k\frac{\log{N}}{\log{y}}+O\left(\frac{1}{\log{y}}\right)-\sum_{j=1}^k\sum_{p|L_j(n)}\left(1-\frac{\log{p}}{\log{y}}\right)\right)+O(S').
\end{align}
Here
\begin{align}
S'=\sum_{\substack{N\le n< 2N\\(\Pi(n),P(z))=1\\\Pi(n)\text{ not square-free}}}\log{N}
\end{align}
We reverse the order of summation over $p$ and $n$, and split the contribution up depending on whether the terms are positive or negative.
\begin{align}
S&=\left(\tau-k\frac{\log{N}}{\log{y}}+O\left(\frac{1}{\log{y}}\right)\right)\sum_{\substack{N\le n< 2N\\ (\Pi(n),P(z))=1}}1-\sum_{j=1}^k\sum_{z\le p< y}\left(1-\frac{\log{p}}{\log{y}}\right)\sum_{\substack{N\le n< 2N\\ (\Pi(n),P(z))=1\\ p|L_j(n)}}1\nonumber\\
&\qquad+\sum_{j=1}^k\sum_{y\le p< 2a_j N+b_j}\left(\frac{\log{p}}{\log{y}}-1\right)\sum_{\substack{N\le n< 2N\\ (\Pi(n),P(z))=1\\ p|L_j(n)}}1+O(S_1).\label{eq:SSplit1}
\end{align}
We notice that the inner sum in the final term makes a contribution only if $L_j(n)$ is a multiple of $p$ and has all prime factors of size at least $z$. Therefore either $p\le L_j(n)/z$ and $L_j(n)$ has prime factors other than $p$, or $L_j(n)=p$. We split the term depending on which of these two is the case.

Thus
\begin{align}
\sum_{y\le p< 2a_jN+b_j}\left(\frac{\log{p}}{\log{y}}-1\right)\sum_{\substack{N\le n< 2N\\ (\Pi(n),P(z))=1\\ p|L_j(n)}}1&=\sum_{y\le p< 2a_jN/z}\left(\frac{\log{p}}{\log{y}}-1\right)\sum_{\substack{N\le n< 2N\\ (\Pi(n),P(z))=1\\ p|L_j(n)}}1+O(1)\nonumber\\
&\qquad+\left(\frac{\log{N}}{\log{y}}-1+O\left(\frac{1}{\log{y}}\right)\right)\sum_{\substack{N\le n< 2N\\ (\Pi(n),P(z))=1\\L_j(n)\text{ prime}}}1.
\end{align}
Substituting this into \eqref{eq:SSplit1} gives
\begin{equation}
S=S_1-S_2+S_3+S_4+O(S')+O(1),
\end{equation}
where
\begin{align}
S_1&=\left(\tau-k\frac{\log{N}}{\log{y}}+O\left(\frac{1}{\log{y}}\right)\right)\sum_{\substack{N\le n< 2N\\ (\Pi(n),P(z))=1}}1,\\
S_2&=\sum_{j=1}^k\sum_{z\le p< y}\left(1-\frac{\log{p}}{\log{y}}\right)\sum_{\substack{N\le n< 2N\\ (\Pi(n),P(z))=1\\ p|L_j(n)}}1,\\
S_3&=\sum_{j=1}^k\sum_{y\le p< 2a_j N/z}\left(\frac{\log{p}}{\log{y}}-1\right)\sum_{\substack{N\le n< 2N\\ (\Pi(n),P(z))=1\\ p|L_j(n)}}1,\\
S_4&=\left(\frac{\log{N}}{\log{y}}-1+O\left(\frac{1}{\log{y}}\right)\right)\sum_{j=1}^k\sum_{\substack{N\le n< 2N\\ (\Pi(n),P(z))=1\\L_j(n)\text{ prime}}}1.
\end{align}
We wish to use the Diamond-Halberstam-Richert $k$-dimensional sieve to get lower bounds for $S_1$, $S_3$ and $S_4$ and an upper bound for $S_2$. We use a simple upper bound to show $S'$ is negligible. This will then give a lower bound for our sum $S$.

We summarize these results in the following proposition.
\begin{prpstn}\label{prpstn:SBounds}
Let $N>N_0$ and $u,v$ be constants satisfying $\tau>u k$, $1<v/(v-1)<u<v$. Let
\[V(z)=\prod_{\substack{p<z\\p\nmid A}}\left(1-\frac{k}{p}\right),\qquad y=N^{1/u},\qquad z=N^{1/v}.\]
Then we have
\begin{align}
S_1&\ge (\tau-k u)NV(z)f_{k}(v)+O\left(\frac{N}{(\log{N})^k\log\log{N}}\right),\\
S_2&\le k N V(z)\int_{1/v}^{1/u}\min\left(F_k(v-vs),\: e^\gamma F_{k-1}\left(\frac{v}{2}\right)w(v-vs)\right)\frac{1-us}{s}d s\nonumber\\
&\qquad+O\left(\frac{N}{(\log{N})^k\log\log{N}}\right),\\
S_3&\ge k N V(z)\int_{1/u}^{1-1/v}\max\left(f_k(v-vs),\: e^\gamma f_{k-1}\left(\frac{v}{2}\right)w(v-vs)\right)\frac{us-1}{s}ds\nonumber\\
&\qquad+O\left(\frac{N}{(\log{N})^k\log\log{N}}\right),\\
S_4&\ge k N V(z)\frac{(u-1)e^\gamma f_{k-1}\left(\frac{v}{2}\right)}{v}+O\left(\frac{N}{(\log{N})^k\log\log{N}}\right),\\
S'&\ll N^{1-1/2v}.
\end{align}
Here $f_k$ and $F_k$ are the lower and upper sieve functions of the Diamond-Halberstam-Richert sieve of dimension $k$, and $w(u)$ is the Buchstab function.
\end{prpstn}
We now establish Proposition \ref{prpstn:MainProp} from Proposition \ref{prpstn:SBounds}.

Given $u,v$ satisfying the conditions of Proposition \ref{prpstn:SBounds} we then have that for $N$ sufficiently large
\begin{equation}
S\ge NV(z)\left((\tau-k u)f_k(v)-k I_1+k I_2+\frac{(u-1)k e^\gamma}{v}f_{k-1}\left(\frac{v}{2}\right)\right)+o\left(\frac{N}{(\log{N})^k}\right),
\end{equation}
where
\begin{align}
I_1&=\int_{1/v}^{1/u}\min\left(F_k\left(v\left(1-s\right)\right),\: e^\gamma F_{k-1}\left(\frac{v}{2}\right)w\left(v\left(1-s\right)\right)\right)\frac{1-us}{s}ds,\\
I_2&=\int_{1/u}^{1-1/v}\max\left(f_k\left(v\left(1-s\right)\right),\: e^\gamma f_{k-1}\left(\frac{v}{2}\right)w\left(v\left(1-s\right)\right)\right)\frac{us-1}{s}ds.
\end{align}
We have that
\begin{equation}
V(z)\gg \frac{N}{(\log{N})^{k}},
\end{equation}
and so
\begin{equation}
S\ge NV(z)\left((\tau-k u)f_k(v)-k I_1+k I_2+\frac{(u-1)k e^\gamma}{v}f_{k-1}\left(\frac{v}{2}\right)+o(1)\right).
\end{equation}
If $v>\beta_k$ then $f_k(v)>0$. Thus we have that $S>0$ provided $N$ is sufficiently large and provided
\begin{equation}
\tau>u k+\frac{k}{f_k(v)}\left(I_1-I_2-\frac{(u-1)e^\gamma}{v}f_{k-1}\left(\frac{v}{2}\right)\right)
\end{equation}
Therefore under these assumptions the $k$-tuple has at most $r_k=\lfloor \tau\rfloor$ prime factors infinitely often. Thus Proposition \ref{prpstn:MainProp} holds.
\section{Proof of Proposition \ref{prpstn:SBounds}}
We first let $N$ be sufficiently large so that $L_i(n)$ is strictly increasing for $n\ge N$ for all $1\le i\le k$ (this happens because $a_i>0$ for all $i$ by Hypothesis \ref{hypthss:Normalised}). In particular $\Pi(n_1)\ne \Pi(n_2)$ for any $n_1\ne n_2$ with $n_1,n_2\ge N$. This assumption is not strictly neccessary, but simplifies notation since we do not have to address set/sequence issues.

We first consider $S_1$. The sum in $S_1$ is already of the correct form to be estimated. We let
\begin{align}
\mathcal{A}&=\{\Pi(n):N\le n< 2N\},\\
\mathcal{A}_d&=\{a\in\mathcal{A}:a\equiv 0\pmod{d}\},
\end{align}
and $\mathcal{P}$ be the set of primes.

We use the standard sieve notation $S(\mathcal{B},\mathcal{Q},z)$ to denote the number of elements of the set $\mathcal{B}$ which are coprime to all the primes in the set $\mathcal{Q}$ that are less than $z$. We see that
\begin{equation}
\sum_{\substack{N\le n< 2N\\ (\Pi(n),P(z))=1}}1=S(\mathcal{A},\mathcal{P},z).
\end{equation}
By virtue of our normalization from Hypothesis \ref{hypthss:Normalised} we have that
\begin{equation}
\#\mathcal{A}_d=g(d)N+O(k^{\omega(d)}),
\end{equation}
where $g(d)$ is the multiplicative function defined by
\begin{equation}
g(p)=\begin{cases}
k,\qquad &p\nmid A,\\
0,&p|A.
\end{cases}
\end{equation}
Thus applying \cite{DiamondHalberstam}[Theorem 9.1] (with the $y$ from their notation taken to be $N(\log{N})^{-6k}$) and recalling that $z=N^{1/v}$ we obtain
\begin{align}
S(\mathcal{A},\mathcal{P},z)&\ge N\prod_{\substack{p<z\\p\nmid A}}\left(1-\frac{k}{p}\right)\left(f_k\left(v-6k\frac{\log\log{N}}{\log{N}}\right)+O\left(\frac{(\log\log{N})^2}{(\log{N})^{1/(2k+2)}}\right)\right)\nonumber\\
&\qquad+O\left(\sum_{m<N(\log{N})^{-6k}}\mu^2(m)(4k)^{\omega(m)}\right),
\end{align}
where $f_k$ is the Diamond-Halberstam-Richert lower sieve function of dimension $k$. 

We see that
\begin{align}
\sum_{m<N(\log{N})^{-2k}}\mu^2(m)(4k)^{\omega(m)}&\le N(\log{N})^{-6k}\sum_{m<N}\frac{\mu^2(m)(4k)^{\omega(m)}}{m}\nonumber\\
&\le N(\log{N})^{-6k}\prod_{p<N}\left(1+\frac{4k}{p}\right)\nonumber\\
&\ll N(\log{N})^{-2k}.
\end{align}
By our construction of $A$ have that
\begin{equation}
V(z)=\prod_{\substack{p<z\\p\nmid A}}\left(1-\frac{k}{p}\right)\asymp (\log{z})^{-k}\asymp (\log{N})^{-k},
\end{equation}
and so the error therm contributes a negligible amount.

By \cite{DiamondHalberstam}[Theorem 6.1] $f_k(x)$ satisfies a Lipschitz condition for $x\ge 1$. Therefore
\begin{equation}
S(\mathcal{A},\mathcal{P},z)\ge NV(z)\left(f_k(v)+O\left(\frac{1}{\log\log{N}}\right)\right).
\end{equation}
We recall that $y=N^{1/u}$ for some fixed constant $u$. Thus provided $\tau\ge k u$ we have 
\begin{align}
S_1&=\left(\tau-ku+O\left(\frac{1}{\log{N}}\right)\right)S(\mathcal{A},\mathcal{P},z)\nonumber\\
& \ge (\tau-k u)NV(z)f_{k}(v)+O\left(\frac{N}{(\log{N})^k\log\log{N}}\right).
\end{align}
We now consider the sum $S_2$. We will obtain two different upper bounds for the inner sum in $S_2$, one of which will give stronger results when $p$ is small, and the other will cover the case when $p$ is large. We first note that
\begin{align}
S_2&=\sum_{z\le p< y}\left(1-\frac{\log{p}}{\log{y}}\right)\sum_{j=1}^k\sum_{\substack{N\le n< 2N\\ (\Pi(n),P(z))=1\\ p|L_j(n)}}1\nonumber\\
&=\sum_{z\le p< y}\left(1-\frac{\log{p}}{\log{y}}\right)S(\mathcal{A}_p,\mathcal{P},z).
\end{align}
We can use \cite{DiamondHalberstam}[Theorem 9.1] to give an upper bound in the same manner as our bound for $S_1$. This gives (using the $y$ from their notation as $N(\log{N})^{-2k-5}/p$)
\begin{align}
S(\mathcal{A}_p,\mathcal{P},z)&\le \frac{k}{p}NV(z)\left(F_k\left(v-v\frac{\log{p}}{\log{N}}\right)+O\left(\frac{1}{\log\log{N}}\right)\right)\nonumber\\
&\qquad+O\left(\sum_{d\le N(\log{N})^{-2k-5}/p}\mu^2(d)4^{\omega(d)}\right).
\end{align}
Thus summing over $p\in[P,2P)$ and treating the error term as before we obtain
\begin{align}
\sum_{P\le p< 2P}\left(1-\frac{\log{p}}{\log{y}}\right)S(\mathcal{A}_p,\mathcal{P},z)&\le \frac{k\log{2}}{\log{P}}\left(1-\frac{\log{P}}{\log{y}}\right)NV(z)F_k\left(v-v\frac{\log{P}}{\log{N}}\right)\nonumber\\
&\qquad+O\left(\frac{N}{(\log{P})(\log{N})^k\log\log{N}}\right).\label{eq:S2Bound1}
\end{align}
When $P$ is small the above bound provides a good estimate, but when $P$ is large we can do better. We define
\begin{align}
\Pi^{(j)}(n)&=\prod_{i\ne j}L_i(n),\\
\mathcal{A}^{(j)}&=\{\Pi^{(j)}(n):N\le n< 2N,(L_j(n),P(z))=1\},\\
\mathcal{A}^{(j,d)}&=\{\Pi^{(j)}(n):N\le n< 2N,(L_j(n),P(z))=1, L_j(n)\equiv 0 \pmod{d}\}.
\end{align}
We then see that, since our forms are coprime by Hypothesis \ref{hypthss:Normalised} we have
\begin{equation}
S(\mathcal{A}_p,\mathcal{P},z)=\sum_{j=1}^kS(\mathcal{A}^{(j,p)},\mathcal{P},z).\label{eq:S2k-1Terms}
\end{equation}
The terms $S(\mathcal{A}^{(j,p)},\mathcal{P},z)$ correspond to $k-1$ dimensional sieves rather than a $k$ dimensional sieve. Reducing the sieve dimension in this way can improve estimates, but the set of $n\in [N,2N)$ such that $(L_j(n),P(z))=1$ does not have as strong level of distribution results, which means that this step is only useful when $p$ is relatively large. 

Since any $L_j(n)$ can have at most $\lfloor v\rfloor$ prime factors if $(L_j(n),P(z))=1$ and $N$ is sufficiently large, we have
\begin{equation}
\sum_{P\le p< 2P}S(\mathcal{A}^{(j,p)},\mathcal{P},z)=\sum_{r=1}^{\lfloor v+1\rfloor}S(\mathcal{B}^{(j)}_{P,r},\mathcal{P},z)+O\left(S'_{P}\right),\label{eq:S2FirstDecomposition}
\end{equation}
where
\begin{align}
\mathcal{B}_{P,r}^{(j)}=&\Bigl\{\Pi^{(j)}(n):N\le n< 2N, (L_j(n),P(z))=1,\nonumber\\
&\qquad L_j(n) \text{ has at least $r$ prime factors in $[P,2P)$}\Bigr\},\\
S'_P&=\sum_{\substack{N\le n< 2N\\ (\Pi(n),P(z))=1\\ p^2|\Pi(n)\text{ for some }p\in[P,2P)}}1.
\end{align}
We define
\begin{equation}
f_{P,r}(n)=\begin{cases}
1,\qquad&\text{$n$ has at least $r$ prime factors in $[P,2P)$ and $(n,P(z))=1$,}\\
0,&\text{otherwise,}
\end{cases}\\
\end{equation}
To estimate $S(\mathcal{B}_{P,r}^{(j)},\mathcal{P},z)$ we wish to estimate the number of elements $\#(\mathcal{B}_{P,r}^{(j)})_d$ of $\mathcal{B}_{P,r}$ which are a multiple of some integer $d$.

We see that $\#(\mathcal{B}_{P,r}^{(j)})_d$ is $0$ unless $(d,A)=1$ by Hypothesis \ref{hypthss:Normalised}. If $(d,A)=1$ then we have
\begin{align}
\#(\mathcal{B}_{P,r}^{(j)})_d&=\sum_{\substack{N\le n< 2N\\ d|\Pi^{(j)}(n)}}f_{P,r}(L_j(n))\\
&=\sum_{\substack{d_1 \dots d_k=d\\ d_j=1}}\sum_{\substack{N\le n< 2N\\ d_i|L_i(n)\:\forall i}}f_{P,r}(a_jn+b_j)
\end{align}
We let $m=a_jn+b_j$ so $a_jN+b_j\le m< 2a_jN+b_j$ and $m\equiv b_j\pmod{a_j}$. The condition $d_i|L_i(n)$ introduces the condition $a_i m\equiv a_ib_j-a_jb_i\pmod{d_i}$ since $(a_j,d_i)=1$ (because $(d,A)=1$). We combine all these congruence conditions via the Chinese Remainder Theorem to give $m\equiv m_0\pmod{a_jd}$ for some $m_0$. By Hypothesis \ref{hypthss:Normalised} we see that $m_0$ is coprime to $a_j d$. Thus
\begin{align}
\#(\mathcal{B}_{P,r}^{(j)})_d&=\sum_{\substack{d_1 \dots d_k=d\\ d_j=1}}\sum_{\substack{a_jN+b_j\le m<2a_jN+b_j\\ m\equiv m_0 \pmod{a_jd}}}f_{P,r}(m)\\
&=\sum_{\substack{d_1 \dots d_k=d\\ d_j=1}}\left(\frac{X_{P,r}^{(j)}}{\phi(d)}+O(E_{P,r}(d))\right).
\end{align}
Here
\begin{align}
X_{P,r}^{(j)}&=\frac{1}{\phi(a_j)}\sum_{a_j N+b_j\le m< 2a_j N+b_j}f_{P,r}(m),\\
E_{P,r}(q)&=\max_{(a,q)=1}\Biggl|\sum_{\substack{a_j N+b_j\le m<2a_j N+b_j\\ m\equiv a \pmod{q}}}f_{P,r}(m)-\frac{\phi(a_j)}{\phi(q)}X_{P,r}^{(j)}\Biggr|.
\end{align}
The function $f_{P,r}$ is a sum of $O(1)$ characteristic functions of numbers with a fixed number of prime factors (at most $\lfloor v\rfloor$) where the prime factors lie in specific intervals (each factor is prescribed to be in one of  $[z,P)$, $[P,2P)$ or $[2P,N)$). By Motohashi \cite{Motohashi} all these characteristic functions have level of distribution equal to $1/2$. Therefore all of the $f_{P,r}$ have level of distribution equal to $1/2$. In particular for any constant $C$ there exists a constant $C'(C)$ such that
\begin{equation}
\sum_{q\le N^{1/2}(\log{N})^{-C'(C)}}\mu^2(q)(4k)^{\omega(q)}E_{P,r}(q)\ll N(\log{N})^{-C}.\label{eq:ErrorBound}
\end{equation}
Thus if
\begin{equation}
(\mathcal{B}_{P,r}^{(j)})_d=\{b\in\mathcal{B}_{P,r}^{(j)}:b\equiv 0\pmod{d}\}
\end{equation}
then from Hypothesis \ref{hypthss:Normalised} we have
\begin{equation}
\#(\mathcal{B}_{P,r}^{(j)})_d=g(d)X_{P,r}^{(j)}+O((k-1)^{\omega(d)}E_{P,r}(a_jd))
\end{equation}
where $g(d)$ is the multiplicative function defined by
\begin{equation}
g(p)=\begin{cases}
\frac{k-1}{p-1},\qquad &p\nmid A,\\
0,&\text{otherwise.}\label{eq:gDef}
\end{cases}
\end{equation}
With this we can now apply \cite{DiamondHalberstam}[Theorem 9.1]. We obtain (taking the $y$ from their notation to be $N^{1/2}(\log{N})^{-C'(2k)}$)
\begin{align}
S(\mathcal{B}_{P,r}^{(j)},\mathcal{P},z)&\le X^{(j)}_{P,r}\prod_{\substack{p<z\\p\nmid A}}\left(1-\frac{k-1}{p-1}\right)\left(F_{k-1}\left(\frac{v}{2}\right)+O\left(\frac{\log\log{N}}{\log{N}}\right)\right)\nonumber\\
&\qquad+O\Bigl(\sum_{q\le N^{1/2}(\log{N})^{-C'(2k)}}\mu^2(q)(4k)^{\omega(q)}E_{P,r}(q)\Bigr).\label{eq:S2UpperBound}
\end{align}
We see from \eqref{eq:ErrorBound} that last term is $O(N(\log{N})^{-2k})$.

We notice that by Merten's theorem
\begin{align}
\prod_{\substack{p<z\\p\nmid A}}\left(1-\frac{k-1}{p-1}\right)&=\prod_{\substack{p<z\\p\nmid A}}\left(1-\frac{k}{p}\right)\left(1-\frac{1}{p}\right)^{-1}\nonumber\\
&=V(z)\left(\frac{e^{\gamma}\phi(A)\log{z}}{A}+O(1)\right).\label{eq:S2VzTerm}
\end{align}
We are therefore left to evaluate the terms $X^{(j)}_{P,r}$. We see that
\begin{align}
\sum_{r=1}^{\lfloor v+1\rfloor}X^{(j)}_{P,r}&=\frac{1}{\phi(a_j)}\sum_{\substack{a_j N+b_j\le m< 2a_j N+b_j\\ (m,P(z))=1}}\sum_{\substack{P\le p< 2P\\p|m}}1\nonumber\\
&=\frac{1}{\phi(a_j)}\sum_{\substack{a_j N\le m< 2a_j N\\ (m,P(z))=1}}\sum_{\substack{P\le p< 2P\\p|m}}1+O(1)\nonumber\\
&=\frac{1}{\phi(a_j)}\sum_{P\le p< 2P}\sum_{\substack{a_j N/p\le n< 2a_j N/p\\ (n,P(z))=1}}1+O(1)\label{eq:S2MainTerm1}
\end{align}
We can evaluate the inner sum asymptotically. By \cite{MontgomeryVaughan}[Theorem 7.11] we have
\begin{equation}
\sum_{\substack{a_j N/p\le n< 2a_j N/p\\ (n,P(z))=1}}1=w\left(v-v\frac{\log{p}}{\log{N}}\right)\frac{a_jN}{p\log{z}}+O\left(\frac{N}{p(\log{N})^2}\right),\label{eq:S2MainTerm2}
\end{equation}
where $w(u)$ is the Buchstab function defined by the delay differential equation
\begin{align}
w(u)&=u^{-1},&\qquad &\text{for $1\le u\le 2$},\nonumber\\
(u w(u))'&=w(u-1),&&\text{for $u>2$.}\label{eq:BuchstabFn}
\end{align}
Putting together \eqref{eq:S2MainTerm1} and \eqref{eq:S2MainTerm2} we obtain
\begin{align}
\sum_{r=1}^{\lfloor v+1\rfloor}X_{P,r}^{(j)}&=\frac{a_j N}{\phi(a_j)\log{z}}\sum_{P\le p< 2P}\frac{1}{p}w\left(v-v\frac{\log{p}}{\log{N}}\right)+O\left(\frac{N}{(\log{P})(\log{N})^2}\right)\nonumber\\
&=\frac{a_jN\log{2}}{\phi(a_j)(\log{P})(\log{z})}w\left(v-v\frac{\log{P}}{\log{N}}\right)+O\left(\frac{N}{(\log{P})(\log{N})^2}\right).\label{eq:S2MainTerm3}
\end{align}
Substituting \eqref{eq:S2VzTerm}, \eqref{eq:S2MainTerm3} and \eqref{eq:S2UpperBound} into \eqref{eq:S2FirstDecomposition} we obtain
\begin{align}
\sum_{P\le p< 2P}S(\mathcal{A}_p^{(j)},\mathcal{P},z)&\le \frac{e^\gamma\phi(A)a_j\log{2}}{A\phi(a_j)\log{P}}NV(z)F_{k-1}\left(\frac{v}{2}\right)w\left(v-v\frac{\log{P}}{\log{N}}\right)\nonumber\\
&\qquad+O\left(\frac{N}{(\log{P})(\log{N})^{k}(\log\log{N})}\right)+O(S'_P).\label{eq:S2MainTerm4}
\end{align}
By Hypothesis \ref{hypthss:Normalised}, $a_j$ and $A$ have the same prime factors (ignoring multiplicity). Thus $\phi(A)a_j=A\phi(a_j)$. Using this, and combining \eqref{eq:S2MainTerm4} with \eqref{eq:S2k-1Terms} we obtain
\begin{align}
\sum_{P\le p< 2P}S(\mathcal{A}_p,\mathcal{P},z)&\le \frac{k e^\gamma\log{2}}{\log{P}} NV(z)F_{k-1}\left(\frac{v}{2}\right)w\left(v-v\frac{\log{P}}{\log{N}}\right)\nonumber\\
&\qquad+O\left(\frac{N}{(\log{P})(\log{N})^{k}(\log\log{N})}\right).\label{eq:S2MainTerm5}
\end{align}
To ease notation we let $s=\log{P}/\log{N}$ and recall that $y=N^{1/u}$. Combining \eqref{eq:S2MainTerm5} with \eqref{eq:S2Bound1} gives
\begin{align}
\sum_{P\le p< 2P}&\left(1-\frac{\log{p}}{\log{y}}\right)S(\mathcal{A}_p,\mathcal{P},z)\nonumber\\
&\le \frac{k\left(1-us\right)NV(z)\log{2}}{s\log{N}}\min\left(F_k(v-vs),\: e^\gamma F_{k-1}\left(\frac{v}{2}\right)w(v-vs)\right)\nonumber\\
&\qquad+O\left(\frac{N}{(\log{P})(\log{N})^k\log\log{N}}\right)+O(S'_P).
\end{align}
An application of partial summation now gives
\begin{align}
S_2&=\sum_{z\le p< y}\left(1-\frac{\log{p}}{\log{y}}\right)S(\mathcal{A}_p,\mathcal{P},z)\nonumber\\
&\le k N V(z)\int_{1/v}^{1/u}\min\left(F_k(v-vs),\: e^\gamma F_{k-1}\left(\frac{v}{2}\right)w(v-vs)\right)\frac{1-us}{s}ds\nonumber\\
&\qquad+O\left(\frac{N}{(\log{N})^k\log\log{N}}\right)+O(S').\label{eq:S2Bound}
\end{align}
Exactly the same argument allows us to bound $S_3$ from below. Using the lower bounds in place of the upper bounds, we obtain
\begin{align}
S_3&=\sum_{y\le p< AN/z}\left(\frac{\log{p}}{\log{y}}-1\right)S(\mathcal{A}_p,\mathcal{P},z)\nonumber\\
&\ge k N V(z)\int_{1/u}^{1-1/v}\max\left(f_k(v-vs),\: e^\gamma f_{k-1}\left(\frac{v}{2}\right)w(v-vs)\right)\frac{us-1}{s}d s\nonumber\\
&\qquad+O\left(\frac{N}{(\log{N})^k\log\log{N}}\right)+O(S').\label{eq:S3Bound}
\end{align}
We now bound $S_4$. We see that
\begin{align}
S_4&=\left(\frac{\log{N}}{\log{y}}-1+O\left(\frac{1}{\log{N}}\right)\right)\sum_{j=1}^k\sum_{\substack{N\le n< 2N\\ (\Pi(n),P(z))=1\\L_j(n)\text{ prime}}}1\\
&=\left(\frac{\log{N}}{\log{y}}-1+O\left(\frac{1}{\log{N}}\right)\right)\sum_{j=1}^k S(\mathcal{C}^{(j)},\mathcal{P},z),
\end{align}
where
\begin{equation}
\mathcal{C}^{(j)}=\{\Pi^{(j)}(n):N\le n< 2N,L_j(n)\text{ prime}\}.
\end{equation}
We let
\begin{align}
\#(\mathcal{C}^{(j)})_d&=\#\{c\in\mathcal{C}^{(j)}:c\equiv 0 \pmod{d}\}\nonumber\\
&=g(d)\#C^{(j)}+(k-1)^{\omega(d)}E(d),
\end{align}
where $g(d)$ is the multiplicative function defined by \eqref{eq:gDef}.

By the Bombieri-Vingradov theorem, for any constant $C$ there is a $C'=C'(C)$ such that
\begin{equation}
\sum_{q\le N^{1/2}(\log{N})^{-C'}}\mu^2(q)(4k)^{\omega(q)}|E(q)|\ll N(\log{N})^{-C}.
\end{equation}
Thus we can apply \cite{DiamondHalberstam}[Theorem 9.1] (with $y$ from their notation as $N^{1/2}(\log{N})^{-C'}$) to give
\begin{align}
S(\mathcal{C}^{(j)},\mathcal{P},z)&\ge \#\mathcal{C}^{(j)}\prod_{\substack{p<z\\p\nmid A}}\left(1-\frac{k-1}{p-1}\right)f_{k-1}\left(\frac{v}{2}\right)+O\left(\frac{N}{(\log{N})^k\log\log{N}}\right)\nonumber\\
&=\left(\frac{N}{\log{N}}+O\left(\frac{N}{(\log{N})^2}\right)\right)V(z)\left(e^\gamma \log{z}+O(1)\right)f_{k-1}\left(\frac{v}{2}\right)\nonumber\\
&\qquad+O\left(\frac{N}{(\log{N})^k\log\log{N}}\right)\nonumber\\
&=\frac{NV(z)e^\gamma f_{k-1}\left(\frac{v}{2}\right)}{v}+O\left(\frac{N}{(\log{N})^k\log\log{N}}\right).
\end{align}
Thus
\begin{equation}
S_4\ge k N V(z)\frac{(u-1)e^\gamma f_{k-1}\left(\frac{v}{2}\right)}{v}+O\left(\frac{N}{(\log{N})^k\log\log{N}}\right).
\end{equation}
Lastly we bound $S'$. We have
\begin{align}
S'&\ll \log{N}\sum_{\substack{N\le n< 2N\\ \Pi(n)\text{ not square-free}}}1\nonumber\\
&\ll \log{N}\sum_{z\le p\ll N^{1/2}}\sum_{\substack{N\le n< 2N\\p^2|\Pi(n)}}1\nonumber\\
&\ll \log{N}\sum_{z\le p\ll N^{1/2}}\left(\frac{k N}{p^2}+O(1)\right)\nonumber\\
&\ll \frac{N\log^2{N}}{z}+N^{1/2}\nonumber\\
&\ll N^{1-1/2v}.
\end{align}
Thus Proposition \ref{prpstn:SBounds} holds.
\section{Acknowledgment}
I would like to thank my supervisor, Prof Heath-Brown for many helpful comments.
\bibliographystyle{acm}
\bibliography{bibliography}
\end{document}